\documentclass[10pt]{article}
\usepackage{times,amsfonts,amsmath,amssymb,exscale,mathrsfs}
\usepackage{pst-all}

\usepackage[spanish]{babel}

\newenvironment{comentario}{\medskip\noindent\footnotesize{\bf Comentario:\ }}{\(\Box\)\medskip}
\title{La cuadratura gaussiana seg\'un Gauss
}
\author{
J.M. Sanz-Serna\footnote{Departamento de Matem\'aticas, Universidad Carlos III de Madrid, Avenida de la Universidad 30, E-28911 Legan\'es (Madrid). Email: jmsanzserna@gmail.com}
}

\date{\today}
\begin{document}
\maketitle

\noindent{\bf Resumen: }Este art\'{\i}culo es una traducci\'on al castellano, resumida y con comentarios, de
la memoria de 1815  en la que Gauss introdujo las reglas de cuadratura que  llevan su nombre. El trabajo de
Gauss  en nada se parece al tratamiento  de la cuadratura gaussiana que dan hoy los libros de texto. La
memoria original es un ejercicio de virtuosismo matem\'atico, basado en un uso magistral de las series, en que
el problema  se reformula como uno de aproximaci\'on funcional que se resuelve con la ayuda de fracciones
continuas.
\bigskip

%\noindent\textbf{Mathematical Subject Classification (2010)} 65C30, 60H05, 16T05

%\noindent\textbf{Keywords} Stochastic

\section{Introducci\'on}

Este art\'{\i}culo es una traducci\'on al castellano, resumida y con comentarios, de la memoria de 1815 \emph{Methodus nova integralium valores per approximationem inveniendi} (M\'e\-todo nuevo de hallar por aproximaci\'on los valores de las integrales) \cite{gauss}, en la que Gauss
introdujo las reglas de cuadratura que hoy llevan su nombre.

Hasta hace algunos a\~{n}os mis conocimientos sobre la cuadratura gaussiana se limitaban m\'as o menos al material com\'un en los libros de texto.\footnote{Creo que la estudi\'e por primera vez en el cl\'asico texto \cite{ik}, del que hay una edici\'on de Dover.} En 2015  prepar\'e una conferencia en la que se comparaban el m\'etodo de Montecarlo y la cuadratura gaussiana. El primero es un algoritmo simple,  f\'acil de describir, cuya ejecuci\'on, consistente en repetir veces y veces algunos pasos muy sencillos, debe hacerse necesariamente en un ordenador. Por el contrario, la cuadratura gaussiana es un algoritmo cuya concepci\'on requiere matem\'aticas nada triviales y, a cambio, puede ejecutarse en el c\'alculo manual por la sencillez de su implementaci\'on. Para mejor preparar tal conferencia le\'{\i} los trabajos originales de Metropolis y sus coautores y de Gauss. Y es entonces cuando descubr\'{\i} que la memoria de Gauss
 en nada se parece al tratamiento  de la cuadratura gaussiana que dan  los libros de texto, tratamiento
 que para la comodidad del lector he reproducido muy resumido en la Secci\'on 1.1 y que, ahora me doy cuenta,
 se debe a Jacobi (1826) \cite{jacobi}. Al contrario, el trabajo original de Gauss es un ejercicio de virtuosismo matem\'atico, un \emph{tour de force} en que el problema se reformula como uno de aproximaci\'on funcional que se resuelve con la ayuda de fracciones continuas de modo similar al usado para aproximar n\'umeros irracionales por racionales. En lenguaje moderno, Gauss usa la transformada de Cauchy, las aproximaciones de Pad\'e, la funci\'on hipergeom\'etrica y otros muchos elementos del  an\'alisis, en un desarrollo cuyo hilo conductor es un manejo maestro de las series. Incluso el teorema de los n\'umeros primos se halla impl\'{\i}cito en el trabajo. Al mismo tiempo, leer la memoria sirve para ver la evoluci\'on de las matem\'aticas en los dos siglos transcurridos: se acep\-taba a comienzos del XIX que un sistema lineal con tantas ecuaciones como inc\'ognitas tendr\'a soluci\'on \'unica, se supon\'{\i}an  todas las funciones  anal\'{\i}ticas, etc. Finalmente Gauss se muestra como un gran analista num\'erico \emph{avant la page} cuando incluye en la memoria tablas de los nodos y pesos de las f\'ormulas con 16 cifras decimales significativas y cierra la misma con un ejemplo num\'erico detallado que \emph{pone ante los ojos} la eficacia de \emph{su} m\'etodo al compararlo con resultados anteriores de otro autor, en este caso Bessel.

\subsection{Cuadratura num\'erica: breve recordatorio}
Antes de pasar adelante, voy a reproducir de modo sucinto la teor\'{\i}a de la cuadratura gaussiana tal y como
se explica  a los estudiantes.\footnote{V\'ease, por ejemplo, el Cap\'{\i}tulo 7 de \cite{ik}.} Esto servir\'a
tanto para presentar el problema con el que tratamos como para apreciar mejor las diferencias entre el
tratamiento original y el  hoy  com\'un.

En la cuadratura\footnote{Por razones hist\'oricas en an\'alisis num\'erico se sigue denominando \emph{cuadratura }al c\'alculo de integra\-les definidas, reserv\'andose el t\'ermino \emph{integraci\'on} para la soluci\'on de ecuaciones diferenciales.} num\'erica se trata de aproximar \(\int_A^B f(x)dx\) por una \emph{regla de cuadra\-tu\-ra} \(\sum_{j=0}^n w_jf(x_j)\); los \(x_j\) se llaman nodos (o abscisas) de la regla y los \(w_j\) son los pesos (o coeficientes).  Si la regla es exacta para las funciones \(f(x) = 1\), \(f(x)=x\), \(f(x)=x^2\), \dots , \(f(x)=x^\ell\) (equivalentemente para cada polinomio de grado \(\leq \ell\)) se dice que su grado de precisi\'on es \(\geq \ell\). Por razones en las que no puedo entrar ahora,\footnote{V\'ease la discusi\'on en las Secciones 5.3 y 5.4.2 de \cite{diez}.} es conveniente elegir los nodos y pesos para lograr que tal grado sea lo mayor posible.

Supongamos primero dados los \(n+1\) nodos, dos a dos distintos. Al imponer grado \(\geq n\) se obtiene un sistema lineal de \(n+1\) ecuaciones para los \(n+1\) pesos y es f\'acil ver que el mismo tiene soluci\'on \'unica. Por otro lado los \(n+1\) valores \(f(x_j)\) permiten construir un \'unico polinomio \(F\) de grado \(\leq n\) que interpola a la funci\'on \(f\), esto es \(F(x_j)=f(x_j)\),
\(j= 0,\dots,n\). Al calcular la integral  \(\int_A^B F(x)dx\) se llega a una expresi\'on de la forma  \(\sum_{j=0}^n w_jf(x_j)\) con los \(w_j\)  independientes de \(f\); por tanto  tomar \(\int_A^B F(x)dx\) como aproximaci\'on a \(\int_A^B f(x)dx\) define una regla de cuadratura, que por su g\'enesis se llama \emph{interpolatoria}. El grado de esa regla es \(\geq n\), pues si \(f\) es un polinomio de grado \(\leq n\) el interpolante \(F\) coincide con la propia \(f\) y por tanto  \(\int_A^B F(x)dx\) con la verdadera integral \(\int_A^B
f(x)dx\). En resumen: dados los nodos, hay una \'unica regla de grado \(\geq n\) y esa regla es la interpolatoria.

?\lq Qu\'e ocurre si los nodos se pueden elegir libremente? Pongamos \(\Omega = \prod_{j=0}^n (x-x_j)\) y, dado un polinomio \(f(x)\) efectuemos la divisi\'on eucl\'{\i}dea \(f(x) = Q(x)\Omega(x)+R(x)\), donde el resto \(R\) es de grado \(\leq n\). Al cuadrar el polinomio \(R\) con la regla intepolatoria basada en los \(x_j\)  obtenemos el valor exacto \(\int_A^B R(x)dx\) (pues el grado de la regla es mayor o igual que el de \(R\)); al cuadrar \(Q\Omega\), la regla proporciona el valor \(0\) (ya que \(Q\Omega\) se anula en cada nodo). Por consiguiente al cuadrar  \(f\)  con la regla se llega a \(\int_A^B R(x)dx\) y esto coincidir\'a con el verdadero valor
 \[
 \int_A^Bf(x)dx=\int_A^B (Q(x)\Omega(x)+R(x)) dx
 \]
 si y solamente si
\[
\int_A^B Q(x)\Omega(x) dx= 0,
\]
 es decir si \(Q\) y \(\Omega\) son mutuamente \emph{ortogonales} respecto del producto escalar de funciones \(\langle g_1,g_2\rangle = \int_A^B g_1(x)g_2(x)dx\).\footnote{ Ver Cap\'{\i}tulo 5 de \cite{ik}.}  Si \(f\) tiene grado \(\leq 2n+1\), el cociente \(Q\) tiene grado \(\leq n\). Resultados sencillos de \'algebra lineal muestra que hay un \'unico polinomio \(\widehat\Omega\) m\'onico (es decir con coeficiente director unidad) de grado \(n+1\) ortogonal a todos los polinomios de grado \(\leq n\);\footnote{Si el intervalo \([A,B]\) coincide con \([-1,1]\), es f\'acil ver, integrando por partes reiteradamente, que \(\widehat\Omega\) es \((1/c_{n+1})P_{n+1}\) donde \[P_{n+1}(x) = \frac{1}{2^{n+1}(n+1)!} \frac{d^{n+1}}{dx^{n+1}} (x^2-1)^{n+1}\] es el polinomio de Le\-gendre y \(c_{n+1}\) su coeficiente director (los \(P_n\) satisfacen la normalizaci\'on de valer 1 en \(x=1\)).   La expresi\'on de los polinomios de  Legendre como combinaci\'on de monomios \(x^j\) se obtiene sin dificultad mediante la recurrencia de tres t\'erminos\[ (n+1)P_{n+1}(x)=(2n+1)xP_{n}(x)-nP_{n-1}(x)\] que da \(P_{n+1}\) conocidos \(P_{n}\) y \(P_{n-1}\) \cite{ch}.

 El caso de  \([A,B]\) general se lleva al de
 \([-1,1]\) cambiando linealmente de variable.
 }
 adem\'as, se prueba que los \(n+1\) ceros de tal polinomio son dos a dos distintos y reales (de hecho pertencen a \((A,B)\)). Tomando pues los ceros de \(\widehat \Omega\) como nodos se obtiene la \'unica regla de cuadratura que alcanza  grado \(\geq 2n+1\).\footnote{De hecho el grado no puede ser superior a \(2n+1\): \(\widehat \Omega\) no puede ser ortogonal a todos los polinomios de grado \(\geq n+1\) pues en particular no es ortogonal a s\'{i} mismo.} Esta es la regla gaussiana.

\subsection{Contenido del presente  art\'{\i}culo}

Como ya dije al principio, este art\'{\i}culo es una traducci\'on al castellano, comentada, de la memoria de
 Gauss. Los comentarios persiguen fundamentalmente facilitar la comprensi\'on al lector. En ning\'un momento ha sido mi objetivo hacer una aportaci\'on a la historia de las matem\'aticas, trabajo que excede  mis competencias. Por tanto, ni pondr\'e la cuadratura gaussiana en el contexto matem\'atico del inicio del siglo XIX, ni seguir\'e su evoluci\'on posterior.

El texto de las cinco secciones siguientes, salvo por tres peque\~{n}as excepciones que se\~{n}alar\'{e},
sintetiza la memoria de Gauss con fidelidad, manteniendo su notaci\'on y terminolog\'{\i}a. Me he apartado del
original en tres puntos: (i) he numerado las ecuaciones para poder referirme a ellas,  (ii) he a\~nadido
puntos o comas al final de las ecuaciones tal y como es usual en los textos matem\'aticos contempor\'aneos y
(iii) he recuadrado algunas de las f\'ormulas y  puesto en \emph{cursiva} ciertas palabras. En Gauss ninguna
ecuaci\'on lleva n\'umero o recuadro. Mis propios comentarios y las explicaciones adicionales que no figuran
en la memoria original van identificados como {\bf  Comentarios} intercalados en el texto y escritos   con
letras de menor tama\~{n}o, o, si son m\'as breves, entre corchetes [ ] o en notas a pie de p\'{a}gina. En
estas explicaciones he usado libremente la notaci\'on y nomenclatura actuales.

\section{Construcci\'on de reglas  interpolatorias}
La memoria consta de 40 p\'aginas y 23 art\'{\i}culos. Comenzaremos nuestro estudio con los
art\'{\i}culos \S 7--\S 12 (p\'aginas 11--21), donde Gauss construye la regla de cuadratura [interpolatoria]
 con nodos arbitrarios [\eqref{eq:formula} m\'as abajo]. Los art\'{\i}culos \S 1--\S 6 (p\'aginas 3--11)
 son paralelos a los \S 7--\S 12  pero referidos al caso particular de reglas con nodos equiespaciados, que Gauss
  nota ya hab\'{\i}an sido desarrolladas por \emph{el sumo} Newton y Cotes [ver p.\ ej.\ \cite{ik}].

Gauss formula as\'{\i} el problema: \emph{determinar \(  \int y\, dx\) entre l\'{\i}mites dados cuando se conocen varios valores de \(y\)}.

\begin{comentario}
Dos observaciones: (i) no se dispon\'{\i}a a\'{u}n de la notaci\'{o}n \(\int_a^b\) para la integral definida (o como dice Gauss para la integral tomada \emph{inter limites datos}), (ii) tampoco se contaba con la notaci\'{o}n \(y = f(x)\) para representar las funciones, debi\'{e}ndose acudir siempre a la per\'{\i}frasis \lq\lq la va\-ria\-ble \(y\) es funci\'{o}n de la va\-ria\-ble \(x\)\rq\rq. Trat\'andose de valores concretos (\lq\lq determinados\rq\rq\ les llama Gauss), hab\'{\i}a que usar \lq\lq \(y_0\) es el valor determinado que toma \(y\) cuando \(x=x_0\)\rq\rq\ en vez de nuestro \(y_0=f(x_0)\).
\end{comentario}

Si la integral \emph{debe ser sumada} desde
\( x = g\)  hasta  \( x = g+\Delta\), Gauss comienza por introducir la variable auxiliar
 \(t =\frac{x-g}{\Delta}\) en t\'erminos de la cual el valor buscado es \( \Delta \int y\, dt\), con la nueva integral tomada  desde \( t=0\) a \(t=1\). Denota por \(  A\), \(  A^\prime\), \(  A^{\prime\prime}\), \( A^{\prime\prime\prime}\), \dots, \(  A^{(n)}\) los \( n+1\) valores dados y por \( a\), \( a^\prime\), \(  a^{\prime\prime}\), \(  a^{\prime\prime\prime}\), \dots, \(  a^{(n)}\)
los correspondientes valores de \( t\).

El primer paso para Gauss es construir
\( Y\), \emph{la funci\'on algebraica de orden\footnote{Hoy dir\'{\i}amos \lq\lq polinomio de grado\rq\rq. Evidentemente, se trata del polinomio de interpolaci\'on escrito en la forma llamada de Lagrange \cite{ik,diez}.} \( n\)}
\begin{eqnarray}\label{eq:Y}
&& \phantom{+} A\: \frac{
(t-a^\prime)(t-a^{\prime\prime})(t-a^{\prime\prime\prime})\cdots (t-a^{(n)})
}%numerator ends
{
(a-a^\prime)(a-a^{\prime\prime})(a-a^{\prime\prime\prime})\cdots (a-a^{(n)})
}%denominator ends
\\%second term begins
&&   + A^{\prime}\: \frac{
(t-a)(t-a^{\prime\prime})(t-a^{\prime\prime\prime})\cdots (t-a^{(n)})
}%numerator ends
{
(a^\prime-a)(a^\prime-a^{\prime\prime})(a^\prime-a^{\prime\prime\prime})\cdots (a^\prime-a^{(n)})
}\nonumber%denominator ends
\\%third term begins
&&   + A^{\prime\prime}\: \frac{
(t-a)(t-a^{\prime})(t-a^{\prime\prime\prime})\cdots (t-a^{(n)})
}%numerator ends
{
(a^{\prime\prime}-a)(a^{\prime\prime}-a^{\prime})(a^{\prime\prime}-a^{\prime\prime\prime})\cdots (a^{\prime\prime}-a^{(n)})
}\nonumber%denominator endsq
\\
&& + {\rm etc.}\nonumber
\\%last term begins
&&   + A^{(n)}\: \frac{
(t-a)(t-a^{\prime})(t-a^{\prime\prime})\cdots (t-a^{(n-1)})
}%numerator ends
{
(a^{(n)}-a)(a^{(n)}-a^{\prime})(a^{(n)}-a^{\prime\prime})\cdots (a^{(n)}-a^{(n-1)})
}.\nonumber%denominator endsq
\end{eqnarray}
Es manifiesto que los valores de esta funci\'on si \(t\) se hace igual a
\( a\), \( a^\prime\), \(  a^{\prime\prime}\), \(  a^{\prime\prime\prime}\), \dots , \(  a^{(n)}\) coinciden con los correspondientes valores de \(y\), de donde deduce que \(Y\) es id\'entica a \(y\) cada vez que \(y\) sea una funci\'on algebraica entera [polinomio] cuyo orden [grado] no exceda de \(n\).

\begin{comentario} En el \S 2, Gauss, al hacer el desarrollo paralelo a \'este en el caso de nodos equiespaciados, hab\'ia ya notado que si dos polinomios de grado \(\leq n\) coinciden para \(n+1\) valores de su variable coinciden id\'enticamente: su diferencia es id\'enticamente nula por ser un polinomio de grado \(\leq n\) con \(n+1\) ceros. El hecho de que un polinomio de grado \(\leq n\) queda determinado por sus valores en \(n+1\) abscisas distintas se usar\'a alguna vez m\'as en lo que sigue.
\end{comentario}

Para hallar \( \int Y\, dt\)  considera sucesivamente las distintas partes de \(  Y\). Designa el producto
\begin{equation}\label{eq:T}
(t-a)(t-a^\prime)(t-a^{\prime\prime})(t-a^{\prime\prime\prime})\cdots (t-a^{(n)})
\end{equation}
por \(T\), y por desarrollo del producto  hace
\begin{equation}\label{eq:Talpha}
T=
t^{n+1}+\alpha t^{n}+\alpha^\prime t^{n-1}+\alpha^{\prime\prime} t^{n-2}+{\rm etc.}+ \alpha^{(n)}.
\end{equation}
El [polinomio en el] numerador de la fracci\'on que aparece multiplicada por \(A\) en el primer t\'ermino de \(Y\) [en \eqref{eq:Y}] ser\'a \(\frac{T}{t-a}\); los numeradores de las partes siguientes ser\'an \(\frac{T}{t-a^\prime}\), \(\frac{T}{t-a^{\prime\prime}}\), \(\frac{T}{t-a^{\prime\prime\prime}}\), etc. En verdad los denominadores no son otra cosa que los valores determinados de los numeradores si respectivamente se pone \(t=a\), \(t=a^\prime\), \(t=a^{\prime\prime}\), \(t=a^{\prime\prime\prime}\), etc. Denota estos denominadores por \(M\), \(M^{\prime}\), \(M^{\prime\prime}\), \(M^{\prime\prime\prime}\), respectivamente, para que sea
\[
Y = \frac{AT}{M(t-a)}+ \frac{A^\prime T}{M^\prime(t-a^\prime)}+\frac{A^{\prime\prime} T}{M^{\prime\prime}(t-a^{\prime\prime})}+{\rm etc.}
+ \frac{A^{(n)} T}{M^{(n)}(t-a^{(n)})}.
\]
Como se hace \(T=0\) para \(t=a\), se tiene:
\begin{eqnarray*}
T &=& t^{n+1}-a^{n+1} + \alpha (t^n - a^n) + \alpha^\prime (t^{n-1}-a^{n-1})+\alpha^{\prime\prime} (t^{n-2}-a^{n-2})\\&&+{\rm etc.}+\alpha^{(n-1)}(t-a).
\end{eqnarray*}
Dividiendo por \(t-a\) ser\'a\footnote{Recordemos que ya se ha notado que   \(T/(t-a)\) es el polinomio de grado \(n\) de coeficiente director unidad con ceros en \(t=a^\prime\), \(t=a^{\prime\prime}\), \dots, \(t=a^{(n)}\). En la f\'ormula que sigue, Gauss expresa este polinomio, no en t\'erminos de sus ceros, sino en t\'erminos de los coeficientes \(\alpha\),  \dots, \(\alpha^{(n-1)}\) de \(T\) en \eqref{eq:Talpha}.
}
\begin{eqnarray}\label{eq:Tt-a}
 \frac{T}{t-a}  & = &  \ t^n            +a t^{n-1}      + aa t^{n-2}           + a^3t^{n-3}                   +{\rm etc.}  +a^n\\
                    &&  \phantom{t^n}+\alpha t^{n-1}    + \alpha at^{n-2}      +\alpha aat^{n-3}              +{\rm etc.}  +\alpha a^{n-1}\nonumber \\
                    &&  \phantom{t^n +a t^{n-1}}        +\alpha^\prime t^{n-2} +\alpha^\prime at^{n-3}        +{\rm etc.}  +\alpha^\prime a^{n-2}\nonumber\\
                    && \phantom{t^n +a t^{n-1} +\alpha^\prime t^{n-2}}         +\alpha^{\prime\prime} t^{n-3}+{\rm etc.}   +\alpha^{\prime\prime}a^{n-3}\nonumber\\
                    && \phantom{t^n +a t^{n-1} +\alpha^\prime t^{n-2}{}+\alpha^{\prime\prime} t^{n-3}}        +{\rm etc.\ etc.}\nonumber\\
                    && \phantom{t^n +a t^{n-1} +\alpha^\prime t^{n-2}{}+\alpha^{\prime\prime} t^{n-3}}\qquad               +\alpha^{(n-1)}.\nonumber
\end{eqnarray}
El valor de esta funci\'on para \(t=a\) resulta
\[
 =(n+1)a^n+n \alpha a^{n-1}+(n-1)\alpha^\prime a^{n-2}+(n-2) \alpha^{\prime\prime}a^{n-3}+{\rm etc.}+\alpha^{(n-1)}.
\]
De aqu\'{\i} Gauss concluye que \emph{\(M\) es el valor de \(\frac{dT}{dt}\) para \(t=a\)}, y
de la misma manera \(M^\prime\), \(M^{\prime\prime}\), \(M^{\prime\prime\prime}\), etc. ser\'an los valores de \(\frac{dT}{dt}\) para \(t=a^\prime\)
, \(t = a^{\prime\prime}\), \(t=a^{\prime\prime\prime}\), etc.

\begin{comentario}
A esta misma conclusi\'on puede llegarse de manera  m\'as f\'acil de varias formas, por ejemplo derivando con respecto a \(t\) en \eqref{eq:T} y evaluando en \(t=a\) o, m\'as elegantemente, observando que el valor del polinomio \(T/(t-a)\)  en \(t=a\) ha de concidir con el l\'{\i}mite al acercarse \(t\) a \(a\) y este l\'{\i}mite  es, por definici\'on, el valor de la derivada de \(T\) en \(t=a\), ya que \(T\) se anula en \(t=a\). Gauss se daba perfecta cuenta de estar matando moscas a ca\~nonazos, porque nota que la conclusi\'on a la que acaba de llegar \emph{consta de otras maneras}. Si recurre a \eqref{eq:Tt-a} es porque a continuaci\'on va a usar esta misma f\'ormula para hallar la integral de \(T/(t-a)\).
\end{comentario}

Halla el valor de la integral \( \int \frac{T\,dt}{t-a} \) desde \(t=0\) hasta \(t=1\) [integrando en \eqref{eq:Tt-a}],
\begin{eqnarray*}
\frac{1}{n+1} + \frac{a}{n}+\frac{aa}{n-1}+\frac{a^3}{n-2}+{\rm etc.}+ a^n\qquad{}  &&\\
  + \frac{\alpha}{n}+\frac{\alpha a}{n-1}+\frac{\alpha a a}{n-2}+{\rm etc.} +\alpha a^{n-1} &&\\
 +\frac{\alpha^\prime}{n-1}+\frac{\alpha^\prime a}{n-2}+{\rm etc.}+ \alpha^\prime a^{n-2}&&\\
 +\frac{\alpha^{\prime\prime}}{n-2}+{\rm etc.}+\alpha^{\prime\prime}a^{n-3}&&\\
  +{\rm etc.etc.}\qquad{}&&\\
  + \alpha^{(n-1)},
\end{eqnarray*}
 expresi\'on cuyos t\'erminos dispone en el orden siguiente
\begin{eqnarray*}
&&\quad a^n +\alpha a^{n-1} + \alpha^\prime a^{n-2}+\alpha^{\prime\prime}a^{n-3}+{\rm etc.}+\alpha^{(n-1)}\\
&&+\frac{1}{2}(a^{n-1}+\alpha a^{n-2}+\alpha^\prime a^{n-3}+{\rm etc.}+\alpha^{(n-2)})\\
&&+\frac{1}{3}(a^{n-2}+\alpha a^{n-3}+\alpha^\prime a^{n-4}+{\rm etc.}+\alpha^{(n-3)})\\
&&+\frac{1}{4}(a^{n-3}+\alpha a^{n-4}+\alpha^\prime a^{n-5}+{\rm etc.}+\alpha^{(n-4)})\\
&& +{\rm etc.}\\
&&+\frac{1}{n-1} (aa+\alpha a+ \alpha^\prime)\\
&&+\frac{1}{n} (a+\alpha)\\
&& +\frac{1}{(n+1)}.
\end{eqnarray*}
 Llegado a este punto, Gauss afirma que  es \emph{manifiesto}\footnote{
  La serie \eqref{eq:series} que se introduce como \emph{deus ex machina} podr\'{\i}a alternativamente haberse presentado de modo natural, como veremos m\'as adelante.
  Es el desarrollo de Laurent de la funci\'on \(\int_0^1 d\tau/ (t-\tau)=\log(t/(t-1))\) anal\'{\i}tica en el dominio resultante de quitar del plano complejo el intervalo \([0,1]\) de la recta real, es decir el intervalo en el que estamos integrando \(Y\).
}
 [?\lq ?]
que a la misma cantidad se llega si en el producto de \( T\) por la serie infinita
\begin{equation}\label{eq:series}
 t^{-1}+\frac{1}{2}t^{-2}+\frac{1}{3}t^{-3}+\frac{1}{4}t^{-4}+{\rm etc.},
 \end{equation}
descartadas las potencias negativas de \( t\) (o dicho m\'as brevemente en la parte del producto que es funci\'on entera [polin\'omica] de \(t\))
 se reemplaza \( t\) por \( a\).
 Hace
 \begin{equation}\label{eq:main}
\boxed{ T
 (t^{-1}+\frac{1}{2}t^{-2}+\frac{1}{3}t^{-3}+\frac{1}{4}t^{-4}+{\rm etc.}) = T^\prime +T^{\prime\prime}}
 \end{equation}
 donde \(T^\prime\) es la funci\'on entera [polinomio] contenida en este producto y \(T^{\prime\prime}\) la otra parte.\footnote{He puesto un recuadro a las f\'ormulas m\'as importantes. El grado de \(T^\prime\), que \emph{no es la derivada de \(T\)}, es \(n\) ya que \(T\) es de grado \(n+1\).} De este modo ser\'a el valor de la integral
 \(\int \frac{Tdt}{t-a}\) desde \(t=0\) hasta \(t=1\) igual al valor de la funci\'on \(T^\prime\) para \(t=a\). Denota los valores determinados de la funci\'on
 \begin{equation}\label{eq:pesos}
 \frac{T^\prime}{\frac{dT}{dt}}
 \end{equation}
 para \(t=a\), \(t=a^\prime\), \(t=a^{\prime\prime}\), \(t=a^{\prime\prime\prime}\), etc. hasta \(t=a^{(n)}\) respectivamente por
 \(R\), \(R^\prime\), \(R^{\prime\prime}\), \(R^{\prime\prime\prime}\), etc. hasta \(R^{(n)}\),  y la integral \(\int Ydt\) de \(t=0\) a \(t=1\) ser\'a
 \begin{equation}\label{eq:formula}
 \boxed{RA+R^\prime A^\prime+R^{\prime\prime}A^{\prime\prime}+{\rm etc.}+R^{(n)}A^{(n)},}
 \end{equation}
lo que multiplicado por  \( \Delta\) exhibir\'a el valor verdadero o aproximado de \( \int y\, dx\) de \(x=g\) a \(x=g+\Delta\).

Ahora que ha obtenido la regla de cuadratura buscada [\eqref{eq:formula}], Gauss replica todo el desa\-rrollo,  trabajando con la variable auxiliar \(  u=2t-1\) en vez de \( t\). [Con la nueva variable auxiliar el intervalo de integraci\'on es sim\'etrico: de \(u=-1\) a \(u=1\).] La funci\'on
\[
 U = (u-b)(u-b^\prime) (u-b^{\prime\prime})\cdots (u-b^{(n)})
\]
 reemplaza a la funci\'on \(T\) en \eqref{eq:T} y la serie\footnote{Desarrollo de Laurent de \((1/2)\log ((u+1)/(u-1))=(1/2)\int_{-1}^1 d\upsilon/(u-\upsilon)\), funci\'on
 anal\'{\i}tica en el complemento  del intervalo real \([-1,1]\).}
\begin{equation}\label{eq:series2}
\varphi = u^{-1}+\frac{1}{3} u^{-3}+\frac{1}{5} u^{-3}+\frac{1}{5} u^{-5}+etc.
\end{equation}
juega el papel antes desempe\~{n}ado por \eqref{eq:series}. En vez de \eqref{eq:main} ahora descompone \(U\varphi\) como \(U^\prime+U^{\prime\prime}\).
Como aplicaci\'on encuentra los pesos de las f\'ormulas de Newton-Cotes [nodos equiespaciados] trabajando tanto con
 \( t\) como con \( u\).

Seguidamente  (\S 11) resuelve un problema algebraico. Dadas  tres funciones enteras [polinomios] \(Z\), \(\zeta\), \(\zeta^\prime\) con coeficientes racionales\footnote{Notar que en \(  \zeta^\prime\) la prima no indica derivada. Gauss siempre usa la notaci\'on de Leibniz (cociente de diferenciales) y no la de Newton  (puntos o primas).  } de una misma indeterminada \(z\), se desea encontrar una funci\'on entera [polinomio] que coincida con \(Z/\zeta\) cuando para \(z\) se toma una cualquiera de las ra\'{\i}ces de la ecuaci\'on \(\zeta^\prime=0\). El \S 12 da en detalle un ejemplo num\'erico del c\'alculo [de tal polinomio].

\begin{comentario} Por supuesto, a la vista de \eqref{eq:pesos}, este problema algebraico es relevante en nuestro contexto en el caso de que \(T\) tenga coeficientes racionales ---y entonces lo mismo les ocurre a \(T^\prime\) y a \(dT/dt\)--- pues las abscisas \(a\), \(a^\prime\), \dots son las ra\'{\i}ces de la ecuaci\'on polin\'omica \(T=0\). En lenguaje moderno sabemos que el cuerpo \({\mathbb Q}(\xi)\) resultante de adjuntar un irracional \(\xi\) al cuerpo racional, que est\'a dado por las expresiones racionales en \(\xi\), coincide, en el caso de que \( \xi\) sea algebraico, con el conjunto de polinomios \({\mathbb Q}[\xi]\).
\end{comentario}
\section{El error}
Los n\'umeros  \S 13 a \S 14 (p\'aginas 22--24) se dedican al an\'alisis del error en la f\'ormula \eqref{eq:formula}.
Se pone
\begin{equation}\label{eq:error}
Ra^m+R^\prime a^{\prime m}+{\rm etc.}+R^{(n)} a^{(n)m}=\frac{1}{m+1} - k^{(m)}
\end{equation}
donde \(  k^{(m)}\) es la diferencia entre la integral \(\int t^mdt\) desde \(t=0\) a \(t=1\) y el valor aproximado.
Si se define
\[
\Theta = kt^{-1}+k^\prime t^{-2}+k^{\prime\prime}t^{-3}+k^{\prime\prime\prime}t^{-4}+{\rm etc.}
\]
o mejor
 (porque \( k\), \( k^\prime\), hasta \( k^{(n)}\) deben anularse\footnote{
Aqu\'{\i} la memoria tiene una incoherencia menor:  en la  definici\'on de \(\Theta\) precedente se deber\'{\i}a haber incluido un t\'ermino \(k_0t^0\) correspondiente al error al cuadrar \(t^0\) y luego haber notado que  el coeficiente \(k_0\) es nulo por la misma raz\'on que lo son
\( k\), \( k^\prime\), \dots, \( k^{(n)}\).
}%
)%
\[
\Theta = k^{(n+1)}t^{-(n+2)}+k^{(n+2)}t^{-(n+3)}+k^{(n+3)}t^{-(n+4)}+{\rm etc.}\: ,
\]
se tendr\'a entonces, desarrollada cada fracci\'on en serie,
\begin{equation}\label{eq:errorsuma}
\frac{R}{t-a}+\frac{R^\prime}{t-a^\prime}+{\rm etc.}+\frac{R^{(n)}}{t-a^{(n)}} = t^{-1}+\frac{1}{2}t^{-2}
+\frac{1}{3}t^{-3}+{\rm etc.} -\Theta.
\end{equation}
\begin{comentario}\label{cota}Se ha desarrollado \(1/(t-a) = t^{-1}+at^{-2}+\cdots\), etc. Vemos ahora claro
 el significado de  la serie \eqref{eq:series}, que Gauss nos present\'o de modo tan artificioso en el desarrollo
  que condujo a \eqref{eq:formula}:  el coeficiente de \(t^{-m-1}\), \(m = 0,1,\dots,\) en esta serie es el
  valor \(1/(m+1)\) de la integral \(\int_0^1 t^{m} dt\), es decir el \emph{momento} de orden \(m\) de la
   medida \(dt\) en el intervalo \([0,1]\). De modo an\'alogo el momento de orden \(m\) de la medida
    (combinaci\'on de masas puntuales o deltas de Dirac) \(R \delta_a+R^\prime \delta_{a^\prime}+\cdots\) asociada
     con la cuadratura \eqref{eq:formula} es el primer miembro de \eqref{eq:error} y \(k^{(m)}\)
     el momento de orden \(m\) de la medida (con signo) asociada con el error. Imponer exactitud para \(1\),
     \(t\), \dots, es, claro, demandar que la cuadratura reproduzca sin error los momentos de ordern 0, 1,
     \dots
     Representado las sucesiones de momentos por las series en potencias de \(1/t\) que los tienen
      como coeficientes (serie que coincide esencialmente con lo que se llama \emph{tranformada Z} o
       tambi\'en \emph{funci\'on generatriz)}, se pasa de \eqref{eq:error} a
\eqref{eq:errorsuma}.

Dada una medida \(\mu\) en la recta real, a la funci\'on \(M(t)=\int_{-\infty}^\infty d\mu(\tau)/(t-\tau)\) se le llama a veces la \emph{transformada de Cauchy de} \(\mu\) (tambi\'en se usan las denominaciones funci\'on de Markov, transformada de Borel, transformada de Stieljes, etc.). Desarrollando  el integrando en potencias de \(1/t\), se ve f\'acilmente que los coeficientes del desarrollo formal en potencias de \(1/t\) de la transformada son los momentos de la medida.  Adem\'as, si la medida tiene soporte compacto, la transformada es de hecho una funci\'on anal\'{\i}tica en el entorno de \(t=\infty\).
Si
el integrando \(f\) es representable mediante una integral de contorno
\(f(\tau) = (1/(2\pi i))\oint_\gamma dt /(\tau-t)\) en el plano complejo, entonces, cambiando el orden de integraci\'on en las variables \(t\) y \(\tau\), se tiene
 \(\int_{-\infty}^\infty f(\tau)d\mu(\tau)= -(1/(2\pi i))\oint_\gamma M(t) f(t) dt \). Esto permite acotar el error de cuadratura para integrandos anal\'{\i}ticos en t\'erminos de cotas
 de la diferencia, evaluada en el contorno de integraci\'on \(\gamma\), entre las transformadas de Cauchy de la verdadera medida y de la regla de cuadratura. En nuestro caso permite acotar el error de cuadratura en t\'erminos de \(\Theta\).
\end{comentario}

Multiplicando por \(T\) ser\'a [vista la igualdad \eqref{eq:main}]
\[
T\left( \frac{R}{t-a}+\frac{R^\prime}{t-a^\prime}+{\rm etc.}+\frac{R^{(n)}}{t-a^{(n)}}\right) =
T^\prime+T^{\prime\prime} -T\Theta.
\]
Gauss nota que la parte primera [primer miembro] de esta ecuaci\'on es una funci\'on entera [polinomio] de \(t\) de orden
[grado] \(n\), cuyos valores determinados  para \(t=a\), \(t=a^\prime\), \(t=a^{\prime\prime}\), etc. respectivamente son\footnote{
Recordemos que ya se observ\'o que \(M\) es el valor del polinomio \(T/(t-a)\) en \(t=a\),  y an\'alogamente para \(a^\prime\), etc.
}
 \(MR\), \(M^\prime R^\prime\),
\(M^{\prime\prime}R^{\prime\prime}\), etc., que son los mismos que los de la funci\'on \(T^\prime\).\footnote{
Ya que \(R\) y \(M\) son los  valores de \(T^\prime/(dT/dt)\) y \(
dT/dt\) en \(a\),  y an\'alogamente para las dem\'as abscisas.
}
Necesariamente esa parte primera debe ser id\'entica a \(T^\prime\), y de ah\'{\i} \(T^{\prime\prime} = T\Theta\).

\begin{comentario}Se tienen por tanto las relaciones fundamentales (que no figuran en Gauss)
\[
\boxed{\frac{R}{t-a}+\frac{R^\prime}{t-a^\prime}+{\rm etc.}+\frac{R^{(n)}}{t-a^{(n)}}= \frac{T^\prime}{T},}\qquad \boxed{\Theta = \frac{T^{\prime\prime}}{T}.} \qquad  (\star)
\]
Las abscisas determinan  \(T\) en \eqref{eq:T} y, a su vez, \(T\) conduce a \(T^\prime\) y \(T^{\prime\prime}\)  a trav\'es de \eqref{eq:main}. El desarrollo de \(T^{\prime\prime}/T\) proporciona los coeficientes de error como Gauss notar\'a acto seguido. La funci\'on racional \(T^\prime/T\) permite recuperar la regla: tiene polos simples en las abscisas y los co\-rres\-pondientes residuos son los pesos.  En \eqref{eq:pesos} reconocemos la bien conocida f\'ormula para calcular los residuos de un cociente en un polo simple.

Por otro lado, consideremos una regla de cuadratura de la forma \eqref{eq:formula} con nodos \(a\), \(a^\prime\), \dots , y pesos \(R\), \(R^\prime\), \dots por el momento arbitrarios. Es evidente que la regla ser\'a exacta para \(t^m\), \(m = 0,1,\dots,n\), si y solo si los \(n+1\) primeros coeficientes del desarrollo de \(R/(t-a)+ R^\prime/(t-a^\prime)+\cdots\) en potencias de \(1/t\) coinciden con los de \eqref{eq:series}, o equivalentemente, tras multiplicar por \(T\), el polinomio
\[T (R/(t-a)+ R^\prime/(t-a^\prime)+\cdots)\]
coincide con el polinomio \(T^\prime\) definido en \eqref{eq:main}. Por tanto, la primera f\'ormula en \((\star)\) caracteriza los pesos de la regla interpolatoria con nodos en los ceros de \(T\).
 \end{comentario}

 Resulta entonces que los coeficientes \(k^{(n+1)}\), \(k^{(n+2)}\), etc. \emph{pueden obtenerse del desarrollo de la fracci\'on \( \frac{T^{\prime\prime}}{T} \)}.
  Encontrados estos la correcci\'on de nuestro valor aproximado de la integral \(\int ydt\)
ser\'a
\[
k^{(n+1)}K^{(n+1)}+ k^{(n+2)}K^{(n+2)}+{\rm etc.}
\]
si la serie en la que se desarrolla \(y\) es
\[
y = K+K^\prime t+K^{\prime\prime}tt+K^{\prime\prime\prime}t^3+{\rm etc.} \:.
\]

\begin{comentario}Gauss  acaba de usar, sin cuestionarse su validez, dos cosas:  (i) que \(y\) puede desarrollarse en serie (de Taylor) en potencias de \(t\) y m\'as a\'un que el desarrollo converge en todo el intervalo \(0\leq t\leq 1\), (ii) que la serie se puede integrar t\'ermino a t\'ermino.
\end{comentario}
\section{La idea principal}

Los n\'umeros \S 15 a \S 16 (p\'aginas 24--26) contienen la idea principal de la memoria.
Cualesquiera que sean los nodos  \( a\), \( a^\prime\), \dots,  Gauss ha obtenido una f\'ormula \eqref{eq:formula} que es exacta para cada polinomio de grado \(\leq n\). Sin embargo  hab\'{\i}a discutido en detalle en el \S 6 c\'omo, para \(n\) par, la regla de Cotes basada en \(n+1\) abscisas equiespaciadas tambi\'en integra exactamente, por simetr\'{\i}a, la funci\'on \(t^{n+1}\).\footnote{ Un caso bien conocido lo proporciona la regla de Simpson (tres nodos) que, aun bas\'andose en interpolaci\'on por una par\'abola cuadr\'atica, cuadra sin error \(t^3\) (y por tanto todos los polinomios c\'ubicos).}
Gauss se plantea entonces c\'omo elegir los nodos para que los coeficientes de error \( k^{(n+1)}\), \( k^{(n+2)}\), \dots (es decir los coeficientes de  \( t^{-n-1}\), \( t^{-n-2}\), \dots, en \( \Theta\)) se anulen sucesivamente.

\begin{comentario}
Es decir hay que lograr que \(\Theta = T^{\prime\prime}/T \) (o equivalentemente \(T^{\prime\prime}\)) tenga en el infinito \(t=\infty\) un cero de orden tan alto como sea posible. (Quiz\'a conviene recordar del cometario que sigue a la f\'ormula \eqref{eq:errorsuma} que los valores de \(\Theta\) en un contorno de integraci\'on suficientemente grande determinan el error de cuadratura cuando el integrando es anal\'{\i}tico.) Para lograr un cero de orden alto, dados los nodos y por tanto \(T\), Gauss tratar\'a de anular sucesivamente los coeficientes de \(1/t\), \(1/t^2\), etc. en \(T^{\prime\prime}\) definido en \eqref{eq:main}.
Como  dispone de \(n+1\) nodos afirma que  podr\'a anular \(n+1\) coeficientes, lo que har\'a la regla exacta para polinomios de grado \(\leq 2n+1\). Da por sentado que imponer \(n+1\) condiciones a los \(n+1\) coeficientes indeterminados de \(T\) en \eqref{eq:Talpha} los determina   de manera \'unica, sin hacer cuesti\'on de la solubilidad del correspondiente sistema lineal de ecuaciones. Tampoco contempla la posibilidad de que, hallado \(T\), careciera de \(n\) ra\'{\i}ces reales dos a dos distintas.
\end{comentario}

Para encontrar la regla de grado m\'as alto posible todo el negocio se reduce para nuestro autor a encontrar, para el valor dado de \(n\), la funci\'on de la forma \eqref{eq:Talpha} que hace que
en el producto [ver \eqref{eq:main}]
\[
 T
 (t^{-1}+\frac{1}{2}t^{-2}+\frac{1}{3}t^{-3}+\frac{1}{4}t^{-4}+{\rm etc.}),
\]
desarrollado en potencias [enteras] de \(t\), los coeficientes  de \(t^{-1}\), \(t^{-2}\), \(t^{-3}\), \dots, \(t^{-(n+1)}\) se anulen.

\begin{comentario}\label{coment}
De
\[
T(t) \int_0^1 \frac{d\tau}{t-\tau} = \int_0^1\frac{T(t)-T(\tau)}{t-\tau}d\tau + \int_0^1 \frac{T(\tau)d\tau}{\tau-\tau}
\]
y \eqref{eq:main} tendremos
\[ T^{\prime}(t) = \int_0^1 \frac{T(t)-T(\tau)}{t-\tau}d\tau, \qquad T^{\prime\prime}(t) = \int_0^1 \frac{T(\tau)d\tau}{t-\tau}.\qquad (\star\star)
\]
Desarrollando
\[
 T^{\prime\prime} = t^{-1} \int_0^1 T(\tau)d\tau+
t^{-2} \int_0^1 \tau T(\tau)d\tau+\cdots
\]
con lo que la aniquilaci\'on de los coeficientes de \(T^{\prime\prime}\) no es otra cosa que la ortogonalidad de  \(T(\tau)\) a \(1\), \(\tau\), \dots
 \end{comentario}

 Alternativamente, si m\'as place, se puede encontrar la funci\'on de la forma \(U = u^{n+1}+\beta u^n+\beta^\prime u^{n-1}+\beta^{\prime\prime}u^{n-2}+{\rm etc.}\), que hace que  su producto por \(u^{-1}+\frac{1}{3} u^{-3}+\frac{1}{5} u^{-5}+\frac{1}{7} u^{-7}+{\rm etc.}\) [ver \eqref{eq:series2}] est\'e libre de las potencias \(u^{-1}\), \(u^{-2}\), \(u^{-3}\), \(u^{-4}\), \dots, \(u^{-(n+1)}\). El modo posterior, nota Gauss, es [por simetr\'{\i}a]] ligeramente m\'as sencillo porque, como se percibe f\'acilmente, los coeficientes de la funci\'on \(U\) que satisface la condici\'on deben anularse alternativamente y poniendo \(\beta = 0\), \(\beta^{\prime\prime} = 0\),
\(\beta^{iv}=0\), etc. se disminuye el trabajo.

 En el ejemplo m\'as sencillo, \(n= 0\), \'unicamente el coeficiente of \( t^{-1}\) en el producto
\((t+\alpha)(t^{-1} + \frac{1}{2}t^{-2}+\frac{1}{3}t^{-3}+{\rm etc.})\) debe anularse. Y como \'aquel sea
\( \frac{1}{2}+\alpha\), tenemos \( \alpha = -\frac{1}{2}\) o \( T = t-\frac{1}{2}\).\footnote{ Es la \emph{regla del punto medio}, basada en interpolaci\'on por una constante, que cuadra sin error cada polinomio de primer grado.} Tras este primer caso, Gauss lleva a cabo los c\'alculos para las reglas con \(n=1\) y \(n=2\), usando tanto la variable auxiliar \(t\) como la \(u\). Sin embargo \emph{este modo, que conduce a c\'alculos cada vez mas molestos, no se lleva m\'as adelante, sino que se avanza hacia la fuente genuina de la soluci\'on.}

\section{Determinaci\'on de los nodos: un camino mejor}

En los n\'umeros \S 17 a \S 21 (p\'aginas 26--36) se proporciona una manera alternativa de hallar \(T\) (o \(U\)) que no precisa de la resoluci\'on de sistemas lineales. Este modo alternativo usa fracciones continuas y Gauss ofrece un breve recordatorio de las mismas.\footnote{Para los lectores que no est\'en familiarizados con las fracciones continuas, \cite{h} es una referencia muy aconsejable. La relaci\'on entre las fracciones continuas y los polinomios ortogonales que se pondr\'a de manifiesto a continuaci\'on se estudia en el Cap\'{\i}tulo III de \cite{ch}.}

\begin{comentario}
Es bien conocido c\'omo las fracciones continuas se emplean para resolver problemas de aproximaci\'on diof\'antica: aproximar n\'umeros irracionales por racionales o tambi\'en racionales por racionales m\'as sencillos. He aqu\'{\i} un ejemplo. Escribiendo \(\pi = 3.141592\dots = 3+0.141592\dots =
3 + 1/(1/0.141592\dots)\), \(1/0.141592\dots = 7.062513\dots \), se obtiene la muy conocida
aproximaci\'on \(\pi\approx 3+1/7=22/7\). El proceso puede reiterarse con \(0.062513\dots = 1/(1/0.062513\dots) =
1/15.996594\dots\), \(\pi \approx 3+(1/(7+1/15))\), etc. para generar una fracci\'on continua cuyo valor es \(\pi\). Cuando se conoce el desa\-rrollo de un n\'umero en fracci\'on continua, las fracciones reducidas proporcionan aproximaciones racionales  que son \'optimas en un sentido que no precisar\'e.

 De modo similar, tratando con funciones en vez de con n\'umeros reales, las fracciones continuas pueden utilizarse para aproximar una funci\'on dada por funciones racionales. Gauss ha reformulado la cuesti\'on que le ocupa como un problema de aproximaci\'on funcional: aproximar \(t^{-1}+\frac{1}{2}t^{-2}+ \cdots\) por una funci\'on racional \(T^\prime/T\) en el sentido de que \(\Theta\), la diferencia de ambas,  tenga un cero de orden lo mayor posible en \(t=\infty\); los ceros de \(T\) (polos de la fracci\'on) son los nodos buscados. Cuando se trabaja con la variable \(u\) hay que aproximar \eqref{eq:series2} por una funci\'on racional \(U^\prime/U\). Es este problema de aproximaci\'on funcional el que va a ser resuelto con ayuda de fracciones continuas.
\end{comentario}

Propuesta la fracci\'on continua
\[
\varphi =\frac{v}{
w +
\frac{v^{\prime}}{
w^{\prime} +\frac{v^{\prime\prime}}{
w^{\prime\prime} +\frac{v^{\prime\prime\prime}}{
w^{\prime\prime\prime}+{\rm etc.}
} % end of fourth quotient
}%end of third quotient
}%end of second quotient
}%end of first quotient
\]
sean formadas dos series de cantidades \( V\), \( V^\prime\), etc. \( W\), \( W^\prime\), etc., por medio de las f\'ormulas
\begin{align}\label{eq:recurrence}
& V = 0&  & W = 1\\
& V^\prime = v&  & W^\prime = wW\nonumber\\
& V^{\prime\prime} = w^\prime V^\prime+ v^\prime V&  & W^{\prime\prime} = w^\prime W^\prime+v^\prime W\nonumber\\
& V^{\prime\prime\prime} = w^{\prime\prime} V^{\prime\prime}+ v^{\prime\prime} V^\prime&  &
W^{\prime\prime\prime} = w^{\prime\prime} W^{\prime\prime}+v^{\prime\prime} W^\prime\nonumber\\
& V^{iv} = w^{\prime\prime\prime} V^{\prime\prime\prime}+ v^{\prime\prime\prime} V^{\prime\prime}&  &
W^{iv} = w^{\prime\prime\prime} W^{\prime\prime\prime}+v^{\prime\prime\prime} W^{\prime\prime}\nonumber
\end{align}
etc. y ser\'a
\begin{eqnarray*}
 \frac{V}{W} & =& 0\\
\frac{V^\prime}{W^\prime} & =&  \frac{v}{w}\\
\frac{V^{\prime\prime}}{W^{\prime\prime}} & =&  \frac{v}{ w+\frac{v^\prime}{w^\prime} }\\
\frac{V^{\prime\prime\prime}}{W^{\prime\prime\prime}} & =&  \frac{v}{ w+\frac{v^\prime}
{w^\prime+\frac{v^{\prime\prime}}{w^{\prime\prime}} }}
\end{eqnarray*}
y as\'{\i} sucesivamente.

De la serie
\[
\frac{v}{WW^\prime}
-\frac{vv^\prime}{W^\prime W^{\prime\prime}}
+\frac{vv^\prime v^{\prime\prime}}{ W^{\prime\prime}W^{\prime\prime\prime}}
-\frac{vv^\prime v^{\prime\prime} v^{\prime\prime\prime}}{ W^{\prime\prime\prime}W^{iv}}
+{\rm etc.}
\]
\noindent el t\'ermino primero es \( = \frac{V^\prime}{W^\prime}\)

\noindent la suma de los dos t\'erminos primeros es \( =\frac{V^{\prime\prime}}{W^{\prime\prime}}\)

\noindent la suma de los tres t\'erminos primeros es \(  =\frac{V^{\prime\prime\prime}}{W^{\prime\prime\prime}}\)

\noindent la suma de los cuatro t\'erminos primeros es \(  =\frac{V^{iv}}{W^{iv}}\)

\noindent y as\'{\i} sucesivamente, de manera que la serie representa la fracci\'on continua \(\varphi\). Del mismo modo se puede representar [mediante una serie]  la diferencia entre \(\varphi\) y cada una de las fracciones aproximantes  \(  \frac{V^\prime}{W^\prime}\), \( \frac{V^{\prime\prime}}{W^{\prime\prime}}\), \( \frac{V^{\prime\prime\prime}}{W^{\prime\prime\prime}}\), etc.

 Para la funci\'on \(\varphi\) en \eqref{eq:series2} que Gauss desea aproximar, la f\'ormula 33 de su memoria  de 1812 sobre la serie hoy llamada hipergeom\'etrica\footnote{Para la funci\'on hipergeom\'etrica puede verse el Cap\'i{\i}tulo XIV de \cite{ww}. } \emph{Disquisitionum generalium circa seriem infinitam \(1+\frac{\alpha \beta}{1 \cdot\gamma}x+\frac{\alpha (\alpha+1)\beta(\beta+1)}{1\cdot 2\cdot \gamma (\gamma+1)}xx+{\rm etc.}\)} suministra el siguiente desarrollo en fracci\'on continua:
 \[
\frac{1}
{% first denominator begins
u-\frac{\frac{1}{3}}
{%second denominator begins
u - \frac{\frac{2\cdot 2}{3\cdot 5}}
{% third begins
u - \frac{\frac{3\cdot 3}{5\cdot 7}}
{% fourth begins
u-\frac{\frac{4\cdot 4}{7\cdot 9}}
{u -{\rm etc.}}
}% fourth ends
}% third ends
}% second denominator ends
}%first denominator ends
,
\]
de donde se tendr\'a
\[
 v= 1,\:  v^\prime = -\frac{1}{3},\: v^{\prime\prime} = -\frac{4}{15}, \: v^{\prime\prime\prime} = -\frac{9}{35}, \: v^{iv} = -\frac{16}{63},\: {\rm etc.}
 \]
\[ w = w^\prime = w^{\prime\prime}  = w^{\prime\prime\prime}= w^{iv}{\rm etc.} = u.\]
De aqu\'{\i} surgen los valores siguientes
\begin{eqnarray*}
&& V = 0, \quad W = 1,\\
&& V^\prime = 1, \quad W^\prime = u,\\
&& V^{\prime\prime} = u, \quad W^{\prime\prime} = uu-\frac{1}{3},\\
&& V^{\prime\prime\prime} = uu-\frac{4}{15}, \quad W^{\prime\prime\prime} = u^3-\frac{3}{5},
\end{eqnarray*}
etc.

\begin{comentario}
Me paro en \(V^{\prime\prime\prime}\), \(W^{\prime\prime\prime}\) para abreviar. La memoria contiene las expresiones hasta \(V^{vii}\), \(W^{vii}\), sin duda porque luego Gauss las necesita en el ejemplo num\'erico.

Los especialistas en polinomios ortogonales reconocer\'an que estos \(W\), \(W^\prime\), \dots\ son los polinomios ortogonales m\'onicos respecto de la medida \(du\) en el intervalo \([-1,1]\) (esto es los polinomios de Legendre  reescalados para hacerlos m\'onicos). La recurrencia \eqref{eq:recurrence}, que permite el c\'alculo de los denominadores de las fracciones reducidas, coincide con la relaci\'on de recurrencia de tres t\'erminos de los polinomios ortogonales.
Los \(V^\prime\), \(V^{\prime\prime}\), etc. son los
correspondientes \emph{polinomios numeradores} o \emph{polinomios asociados} \cite{ch}.
\end{comentario}

Un poco de atenci\'on aclara que \(V\), \(V^\prime\), \( V^{\prime\prime}\), \( V^{\prime\prime\prime}\), etc., \(W\), \(W^\prime\), \( W^{\prime\prime}\), \(W^{\prime\prime\prime}\), etc. son funciones enteras [polinomios] de la indeterminada \(u\). El t\'ermino alt\'{\i}simo de \(V^{(m)}\) se har\'a \(u^{m-1}\) y
t\'ermino alt\'{\i}simo de \(W^{(m)}\) se har\'a \(u^{m}\).
Si \(\varphi -\frac{V^{(m)}}{W^{(m)}}\)
es convertida en serie descendente [por el procedimiento descrito m\'as arriba], su primer termin\'o ser\'a
\[
\frac{2\cdot 2\cdot 3\cdot 3\cdots m\cdot m\: u^{-(2m+1)}}{3\cdot 3\cdots (2m-1)(2m-1)}.
\]

\begin{comentario}
Consideremos \(\varphi -\frac{P}{Q}\) con \(P\) arbitrario de grado \(m-1\) y \(Q\) arbitrario de grado \(m\).
Hay \(2m+1\) coeficientes intedeterminados en \(P\) y \(Q\), pero s\'olo \(m\) son esenciales porque se puede multiplicar \(P\) y  \(Q\) por cualquier constante no nula. Parece entonces
razonable poder elegir \(P\), \(Q\) para anular los \(2m\) coeficientes de \(u^{-1}\), \dots, \(u^{-2m}\) en esa diferencia. Los polinomios \(V^{(m)}\) \(W^{(m)}\) obtenidos de la fracci\'on continua logran tal anulaci\'on. En lenguaje moderno, \(\frac{V^{(m)}}{W^{(m)}}\) es el aproximante de Pad\'e \((m-1,m)\) de \(\varphi\) (en \(t=\infty\)). Ver \cite{jones} sobre las relaciones entre la aproximaci\'on de Pad\'e y las fracciones continuas.
\end{comentario}

El producto \(\varphi W^{(m)}\) estar\'a compuesto de la funci\'on entera [polinomio] \(V^{(m)}\) y de una serie infinita, cuyo primer t\'ermino es
\[
\frac{2\cdot 2\cdot 3\cdot 3\cdots m\cdot m\: u^{-(m+1)}}{3\cdot 3\cdots (2m-1)(2m-1)}.
\]
As\'{\i} se ha encontrado la funci\'on \(U\) de orden [grado] \(n+1\) que satisface la condici\'on establecida anteriormente de que el producto \(\varphi U\)
est\'e libre de las potencias \(u^{-1}\), \(u^{-2}\), \(u^{-3}\), \dots, \(u^{-(n+1)}\). No ser\'a otra que \(W^{(n+1)}\) y al mismo tiempo se manifiesta que
\(U^\prime\) ser\'a igual a \(V^{(n+1)}\).\footnote{Una conocida relaci\'on entre los polinomios ortogonales y los polinomios asociados/numeradores permite entonces escribir \(U^{\prime}(u) = \int_0^1 (U(u)-U(\upsilon))d\upsilon/(u-\upsilon)\), ver f\'ormula (4.6) en \cite{ch}. La relaci\'on que se corresponde con esta cuando se trabaja con \(t\)--\(T\) en lugar de \(u\)--\(U\) se dedujo en la f\'ormula \((\star\star)\) dentro de un comentario  anterior.}
 \emph{Por tanto para [los nodos] \(b\), \(b^\prime\), \(b^{\prime\prime}\), \dots ,\(b^{(n)}\) tomaremos las ra\'{\i}ces de la ecuaci\'on
\(W^{(n+1)}=0\) y calcularemos los valores de los coeficientes [o pesos] de la manera mostrada arriba y nuestra f\'ormula disfrutar\'a de un orden ascendente a \(2n+1\).}

Para concluir esta parte Gauss da expresiones en forma cerrada para los \(W^{(n)}\),\footnote{Esto es, para los polin\'omios de Legendre m\'onicos.}  discute su expresi\'on en t\'erminos de la funci\'on hipergeom\'etrica y repite todo el desa\-rrollo  usando la variable auxiliar \(t\) en lugar de la \(u\).

\section{Uso de las reglas}

Los dos \'ultimos art\'{\i}culos \S 22-\S 23 (p\'aginas 36–-40) se dedican al uso de las reglas.

Para \(  n = 0,\dots, 6\) (uno a siete nodos),\footnote{Para seis y siete nodos la determinaci\'on de los nodos requiere resolver ecuaciones de grado seis o siete, pero por la simetr\'{\i}a \(u\mapsto -u\) ambas se reducen a ecuaciones c\'ubicas.}
Gauss lista:

\begin{enumerate}
\item Los polinomios \( U\), \( U^\prime\), \(T\), \(T^\prime\).
\item Las abscisas \( a\), \( a^\prime\), \dots\ con 16 cifras decimales significativas.
\item Los pesos \( R\), \( R^\prime\), \dots\ con 16 cifras decimales significativas. (Para \( n\geq 3\) tambi\'en
los correspondientes logaritmos decimales con 10 cifras significativas.)
\item El polinomio con coeficientes racionales que evaluado en los nodos da los pesos.
\item El t\'ermino dominante del desarrollo en serie del error.
\end{enumerate}

\begin{comentario}Por supuesto, los logaritmos de los pesos se usaban para efectuar las productos de \'estos por los valores conocidos de \(y\) que figuran en la regla \eqref{eq:formula}.

M\'as precisamente, la memoria da los valores del logaritmo decimal  de \(10^9 R\), \(10^9 R^\prime\), etc. El artificio de multiplicar cada n\'umero por \(10^9\) presumiblemente servir\'{\i}a para evitar usar logaritmos negativos. (Los logaritmos negativos dificultan el c\'alculo: para multiplicar un n\'umero con logaritmo positivo por un segundo con logaritmo negativo hay que efectuar una \emph{resta} y no una suma. En el bachillerato de mi adolescencia las restas se evitaban hallando complementos, por ejemplo
\(\log_{10} 0.2 = -0.698970\dots\) se expresaba  como \(\log_{10} (1/10)+ \log_{10} 2 = -1+0.301030\dots\), cosa que escrib\'{i}amos \(\bar{1}.301030\dots\))
\end{comentario}

 Finalmente \emph{ponemos ante los ojos la eficacia de nuestro m\'etodo calculando los valores de la integral}
 \(\int\frac{dx}{\log x}\)   desde \( x = 100000\) hasta \( x = 200000\)  por medio de las reglas con uno a siete nodos [subrayo las cifras que no se modifican al incrementar \(n\)]:
\begin{eqnarray*}
&&\underline{ 8}390.394608\\
&&\underline{ 840}5.954599\\
&&\underline{ 8406.2}36775\\
&&\underline{ 8406.242}970\\
&&\underline{ 8406.2431}17\\
&&\underline  {8406.243121}\\
&&\underline  {8406.2431211}
\end{eqnarray*}
Gauss nos informa de que Bessel hab\'{\i}a obtenido, no dice c\'omo, el valor 8406,24312.

\begin{comentario}Mostr\'e una versi\'on preliminar de este trabajo al prof. Ll. N. Trefethen  y comprob\'o con su software CHEBFUN los c\'alculos de Gauss. Gauss lista los productos \(\Delta RA\), \(\Delta R^\prime A^\prime\), etc. con siete cifras decimales tras la coma. Para la regla de siete
 nodos, Trefethen encontr\'o que los valores de esos productos proporcionados por Gauss son t\'{\i}picamente
err\'oneos en una unidad en la cifra menos significativa (esto es hay  error de una diezmillon\'esima aproximadamente).

Como vemos de los subrayados, los errores disminuyen exponencialmente al incrementar el n\'umero de nodos. Es f\'acil probar que esto ocurre para cualquier integrando anal\'{\i}tico.
El intervalo de integraci\'on es sumamente largo, pero, a cambio, el integrando apenas var\'{\i}a  (\((d/dx) \log x = 1/x \leq 10^{-5}\)), lo que explica la bondad de los resultados incluso con muy pocos nodos.

 Por \'ultimo no es aventurado suponer que esta integral era de inter\'es para Gauss en relaci\'on con su conjetura de que \(1/\log x \) es la densidad de los n\'umeros primos (actual teorema de los n\'umeros primos). Por curiosidad: hay 8392 primos en el intervalo considerado.
\end{comentario}

\section*{Agradecimientos}
Antonio Garc\'{\i}a y  Nick Trefethen  comentaron  versiones preliminares del texto. Jorge Arves\'u y Paco Marcell\'an me informaron sobre la transformada de Cauchy y los polinomios asociados respectivamente. Gracias a todos.
 %Este trabajo ha recibido apoyo del proyecto MTM2016-77660-P(AEI/FEDER, UE) funded by MINECO.

%
\end{document}